\documentclass{elsarticle}
\usepackage{amssymb,amsmath,latexsym}
\font\tenscr=rsfs10  scaled \magstep0 \font\sevenscr=rsfs7 scaled \magstep0
\font\fivescr=rsfs5 scaled \magstep0 \skewchar\tenscr='177 \skewchar\sevenscr='177
\skewchar\fivescr='177
\newfam\scrfam \textfont\scrfam=\tenscr \scriptfont\scrfam=\sevenscr
\scriptscriptfont\scrfam=\fivescr
\def\mathscr{\fam\scrfam}

\font\stenscr=rsfs5  scaled \magstep3
\newfam\sscrfam \textfont\sscrfam=\stenscr

\newtheorem{thm}{Theorem}
\newtheorem{prop}{Proposition}

\newtheorem{lemma}{Lemma}
\newtheorem{cor}{Corollary}

\newtheorem{rem}{Remark}
\def\dj{d\kern-0.4em\char"16\kern-0.1em}
\def\Dj{\mbox{\raise0.3ex\hbox{-}\kern-0.4em D}}
\def\ls{\leqslant}
\def\gs{\geqslant}
\def\O{\Omega}
\def\proof{%
  \par\topsep6pt plus6pt
  \trivlist
  \item[\hskip\labelsep\it Proof.]\ignorespaces}
\def\endproof{\qed\endtrivlist}
\expandafter\let\csname endproof*\endcsname=\endproof

\def\qedsymbol{\ifmmode\bgroup\else$\bgroup\aftergroup$\fi
  \vcenter{\hrule\hbox{\vrule height.6em\kern.6em\vrule}\hrule}\egroup}
\def\qed{\ifmmode\else\unskip\nobreak\fi\quad\qedsymbol}

\newcommand{\eql}[2]{\begin{equation} \label{#1} #2 \end{equation}}   
\newcommand{\ld}{\backslash} 
\newcommand{\rd}{/} 

\begin{document}

\journal{}

\title{\Large\bf Solving linear equations in fuzzy quasigroups\tnoteref{t1}}
\tnotetext[t1]{Research of the first author is supported by the Serbian Ministry of Education, Science and Technological Development, Grants ON 174008 and ON 174026 and partially through the joint project 'Algebraic and combinatorial structures with applications' of Serbian and Macedonian Academies of Sciences and Arts. The second and the third authors are supported by the Serbian Ministry of Education, Science and Technological Development,  Grant No.\ 174013.}
\author[sasa]{Aleksandar Krape\v z}
\ead{sasa@mi.sanu.ac.rs}

\author[fsuns]{Branimir \v Se\v selja}
\ead{seselja@dmi.uns.ac.rs}

\author[fsuns]{Andreja Tepav\v cevi\'c}
\ead{andreja@dmi.uns.ac.rs}


\address[sasa]{Mathematical Institute of the Serbian Academy of Sciences and Arts, Belgrade, Serbia}
\address[fsuns]{Department of Mathematics and Informatics, Faculty of Sciences, University of Novi Sad, Serbia}

\begin{abstract}\small
We deal with solutions of classical linear equations $a\cdot x =b$ and $y\cdot a=b$, applying a particular lattice valued fuzzy technique. Our framework is a structure with a binary operation $\,\cdot\,$ (a groupoid), equipped with a fuzzy equality. We call it a fuzzy quasigroup if the above equations have unique solutons with respect to the fuzzy equality. We prove that a fuzzy quasigroup can equivalently be characterized as a structure whose quotients of cut-substructures with respect to cuts of the fuzzy equality are classical quasigroups. Analyzing two approaches to  quasigroups in a fuzzy framework, we prove their equivalence. In addition, we prove that a fuzzy loop (quasigroup with a unit element) which is a fuzzy semigroup is a fuzzy group and vice versa.  Finally, using properties of these fuzzy quasigroups, we give answers to existence of solutions of the mentioned linear equations with respect to a fuzzy equality, and we describe solving procedures.
\end{abstract}

\maketitle


\section{Introduction}
The aim of our paper is to combine a classical mathematical (algebraic) topic with fuzzy structures and techniques in order to deal with solutions of basic linear equations with one binary operation.  Namely, we investigate equations of the form $a\cdot x=b$ and $y\cdot a=b$ in the most general structure with a binary operation $\,\cdot\,$. Classical algebraic structures in which such equations have unique solutions are quasigroups. Needless to say, linear equations appearing in specific problems and in applications are not necessarily situated in a quasigroup structure. Therefore, solutions may not exist, or may not be unique, or the equality of objects may be fuzzy, preventing standard solution procedures.
Our framework is an ordinary structure with a binary operation -- a groupoid,  which a priori does not satisfy any condition (identities, special elements, existence of (unique) solutions of equations...). We equip it with a particular fuzzy (lattice valued) equality which we use instead of the classical equality "$\,=\,$". In this framework we investigate existence of solutions and solutions themselves of the above linear equations.

As mentioned, we use classical groupoids and quasigroups,  for which there is a huge literature in algebra and combinatorics (finite quasigroups are Latin Squares), see e.g., \cite{VD,HOP,Smith}.

Concerning the fuzzy approach, we deal with structures with a specific fuzzy equality.

Our basic tool are $L$-sets, introduced 1979 by Fourman and Scott (\cite{FS}) under the name of $\O$-sets. Intention of the authors was to model intuitionistic logic. An $\O$-set is a nonempty set $A$ equipped with an $\O$-valued equality $E$, where $\O$ is a complete Heyting algebra and $E$ is a symmetric and transitive map from $A^2$ to $\O$. This notion has been further applied to non-classical predicate logics, and also to foundations of Fuzzy Set Theory (\cite{gott,hoh1}).

In our approach $L$ is a complete lattice without additional operations. On the one hand, a complete lattice is not sufficiently rich as a truth values structure in the corresponding fuzzy logic; on the other hand, our research is mostly algebraic, and a complete lattice allows main algebraic notions and properties to be preserved under fuzzification by means of cut sets  ("cutworthy approach", see \cite{KY}).
This approach is widely used for dealing with algebraic topics  (see e.g., \cite{DiNG}, then also \cite{BSAT7,BSAT1}),  and with the lattice-valued topology (starting with \cite{hoso} and many others).  In the recent decades a complete lattice is often replaced by a complete residuated lattice (\cite{Bel1}). A detailed approach to cutworthiness for a particular residuated lattice defined on the unit interval has been presented by B\v elohl\'{a}vek in \cite{Belcut}.

A lattice-valued equality  generalizing the classical one  has been introduced in fuzzy mathematics by H\"{o}hle in \cite{hoh}, and then it was used in investigations of fuzzy functions and fuzzy algebraic structures by many authors, in particular by Demirci (\cite{Dem2}), B\v elohl\'{a}vek and Vychodil (\cite{Bel2}) and others.

Identities were analyzed for $L$-algebras in \cite{BSATie}, and then this approach has been developed in \cite{fi1,fi2,fi3,fugr}.

Quasigroups were investigated in the fuzzy framework in the classical way, as suitable fuzzy subsets of a quasigroup, compatible with the operation, mostly by Dudek, Akram or both (\cite{Dud1,Akr,AkrDud1,AkrDud2}), then by
Alshehri (\cite{Als}), Rosenberg (\cite{Ros}). Our approach is different -- our basic structure is an arbitrary groupoid, not a quasigroup; quasigroups appear on quotient structures over cuts, with respect to cuts of the fuzzy equality.

Concerning linear and other equations in the fuzzy framework, investigations have mostly been (and still are) oriented toward equations over fuzzy sets, fuzzy numbers,  etc. Let us mention the papers \cite{San} by Sanchez and \cite{Buc} by Buckley from the early period, then also the recent paper \cite{MMB} by Mazarbhuiya,  Mahanta and Baruah, dealing with a single binary operation. As mentioned, our topic are classical linear equations while the solving methods are fuzzy.\\

The paper is organized as follows. Preliminary section contains basic definitions concerning classical quasigroups. Next we present the  framework of fuzzy (lattice valued) structures with a particular fuzzy equality. We also list previous results relevant to our investigation. In section Results we first introduce our basic structure, $L$-groupoid over which we investigate solutions of linear equations. If the solutions are unique with respect to a fuzzy equality, we obtain an $L$-quasigroup. We characterize it by quotient structures over cuts, which are ordinary quasigroups. Next we introduce an $L$-equasigroup, which is a structure with a fuzzy equality and, equivalently as in the classical algebra, with three binary operations, fulfilling particular identities. We prove that, with respect to the first operation,  it is an $L$-quasigroup. Using the Axiom of Choice, we also prove the converse, that an $L$-quasigroup can be equipped with two additional operations, so that the new $L$-algebra is an $L$-equasigroup. By a suitable example, we show how classical linear equations can be solved uniquely, up to the fuzzy equality. To complete our investigation, we prove that a fuzzy loop (a quasigroup with a unit element) which is a semigroup is an $L$-group and vice versa.
As an application, we show how our procedure can be applied for solving linear equations in the most general situation, having an arbitrary binary operation and a fuzzy equality arising from the concrete real conditions.

\section{Preliminaries}\label{pr}

\subsection{Quasigroups}
\label{Q}
An \textbf{algebra} is a pair $(A,F)$, where $A$ is a nonempty set and $F$ is a collection of operations on $A$. Here we deal mostly with \textbf{groupoids} - algebras with a single binary operation. In addition, we consider algebras with several binary operations.

There are two standard ways to define quasigroups. One is to consider them as special groupoids:

A groupoid $(Q,\,\cdot\,)$ is a \textbf{quasigroup} if for all $a,b \in Q$\/, both linear equations:
\eql{solv}{a \cdot x = b \qquad\qquad\qquad\qquad y \cdot a = b}
are uniquely solvable for $x,y$\/.

 The other way is to define quasigroups as algebras with three binary operations $\cdot,\ld,\rd$ (called \textbf{multiplication, left division} and \textbf{right division} respectively):

An \textbf{equasigroup} is
an algebra $(Q, \,\cdot\,, \,\ld\,,\, \rd\,)$ which satisfies the following identities:

$Q1:\; y \approx x \cdot (x \backslash y) ;$

$Q2:\; y \approx x \backslash (x \cdot y) ;$

$Q3:\; y \approx (y / x) \cdot x ;$

$Q4:\; y\approx (y \cdot x) / x .$
\begin{thm}\sl\label{equiv} If $(Q,\,\cdot\,)$ is a quasigroup, then $(Q, \,\cdot\,, \,\ld\,,\, \rd\,)$ is an equasigroup, where the additional binary operations $\ld$ and $\rd$ are defined by:
\begin{equation}
a\ld b=c\;\mbox{ iff }\; b=a\cdot c\;\qquad\;\mbox{ and } \;\qquad\;
a\rd b=c \mbox{ iff }\; a=c\cdot b.
\end{equation}

\end{thm}

However, it is important to note that these two kinds of quasigroups have different properties. For example, subquasigroups and homomorphic images of quasigroups need not be quasigroups themselves, while with equasigroups this is always the case.

A quasigroup $(Q,\,\cdot\,)$ with an identity element $e$ is a \textbf{loop}: for every $x\in Q$, $e\cdot x=x\cdot e=x$. For our purposes here, we consider a loop as a structure $(Q,\,\cdot\,, e)$ with the nullary operation in the language, corresponding to the identity element.   Alternatively, an equasigroup is an {\bf eloop} if for all $x,y$, $x\backslash x=y/y$; in this approach $x\backslash x$ serves as the identity element. Finally, a \textbf{group} is an associative loop. A group is often defined as an algebra $(G,\cdot,{\;}^{-1}, e)$ with a binary operation $\,\cdot\,$, unary ${\;}^{-1}$ and a constant $e$, such that the binary operation is associative, $e$ is the identity element, and for every $x\in G$,  $x\cdot x^{-1}=x^{-1}\cdot x=e$.

Basic facts about quasigroups can be found in  e.g., \cite{VD, HOP, Smith}. \\

As usual in algebra, we denote the \textbf{quotient structure}  of an algebra $\mathcal{A}$ over the congruence $\theta$  by $\mathcal{A}{/}\theta$
%
 This denotation is usual and commonly accepted. At the same time, an equasigroup possesses a binary operation denoted by the same symbol $/\;$.
  Still, due to the context, no misunderstanding should arise.\\

We use the following version of the Axiom of Choice:\\

(AC) \textsl{For a collection $\mathcal{X}$ of nonempty subsets of a set $M$, there exists a function $f:\mathcal{X}\rightarrow M$, such that for every $A\in\mathcal{X}$, $f(A)\in A$.}

\subsection{$L$-valued functions and relations}

Throughout the paper, $(L ,\wedge ,\vee , \ls ) $ is a complete lattice with the top and the bottom elements 1 and 0 respectively. It is considered to be the co-domain of all membership functions.

An \textbf{$L$-valued function} $\mu$ on a nonempty set $Q$  is a mapping $\mu:Q\rightarrow L$. It is also called a \textbf{fuzzy set on $Q$}, or a \textbf{fuzzy subset of $Q$}, in particular when the codomain lattice $L$ is known from the context (most often when it is a unit interval $[0,1]$, with respect to the classical order $\,\ls\,$).

For $p\in L$, a \textbf{cut set} or a
$p$-\textbf{cut}  of an $L$-valued function $\mu:Q\rightarrow L$ is a
subset $\mu_p$ of $Q$ which is the inverse image of the principal
filter in $L$, generated by $p$:
$$\mu_p=\mu^{-1}({\uparrow}p)=\{x\in Q\mid\mu (x)\gs p\}.$$

An \textbf{$L$-valued} (binary) \textbf{relation $R$ on $Q$} is an $L$-valued function on $Q^2$, i.e., it is a mapping $R: Q^2\rightarrow L$.

Observe that for $p\in L$, by $R_p$ we denote the cut set for an $L$-valued relation $R$ on $Q$, as defined above:
\[ R_p=\mu^{-1}({\uparrow}p)=\{ (x,y)\in Q^2\mid R(x,y)\gs p.\}\]
\[ R_p=R^{-1}({\uparrow}p)=\{ (x,y)\in Q^2\mid R(x,y)\gs p.\}\]
\noindent $R$ is \textbf{symmetric} if
\begin{equation}R (x,y)= R (y,x)\; \mbox{ for all }\; x,y\in Q \hspace{1.4cm}\label{sym}\end{equation}
%
and \textbf{transitive} if
\begin{equation} R (x,y)\gs  R (x,z)\wedge R (z,y)\mbox{ for all }\; x,y,z\in Q.\label{tran}\end{equation}

Let $\mu:Q\rightarrow L$ and $R: Q^2\rightarrow L$ be an $L$-valued function and an $L$-valued relation on $Q$, respectively. Then we say that $R$ is an \textbf{$L$-valued relation on $\mu$} if for all $x,y\in Q$
\begin{equation} R(x,y)\ls\mu (x)\wedge\mu (y)\label{ronm}.\end{equation}

An $L$-valued relation $ R$ on  $\mu:Q\rightarrow L$ is said to be \textbf{reflexive on $\mu$}, or \textbf{$\mu$-reflexive} if
\begin{equation}\mbox{$ R (x,x)=\mu (x)$ for every $x\in Q$.}\label{rf}\end{equation}

A symmetric and transitive $L$-valued relation $ R$ on $Q$, which is reflexive on $\mu:Q\rightarrow L$ is an \textbf{$L$-valued equivalence} on $\mu$.

An $L$-valued equivalence $R$ on $Q$ fulfills the \textbf{strictness} property (see \cite{hoh1}):
\begin{equation} R (x,y)\ls R(x,x)\wedge R(y,y).\label{str}\end{equation}

More precisely, the strictness property follows from symmetry and transitivity only. The proof is straightforward.

An  $L$-valued equivalence $R$ on $Q$  is an \textbf{$L$-valued equality}, if it satisfies:
\begin{equation}R(x,y)=1 \; \mbox{ implies }\; x=y.\label{nv1}\end{equation}

\begin{rem}\label{refl1}\rm
The above properties of $L$-valued relations are not uniquely defined in the literature.
Firstly, reflexivity as defined here, or $\mu$-reflexivity,  is different from the classical condition $R(x,x)=1$ for all $x\in Q$. The main reason is that  $R:Q^2\rightarrow L$ is considered here to be an $L$-valued relation on a function, i.e., on a fuzzy set $\mu :Q\rightarrow L$. Such $L$-valued relations are supposed to fulfill the property (\ref{ronm}). Therefore the value $R(x,x)$ could not be greater than  $\mu (x)$.
Next, an $L$-valued equality is defined here as an $L$-valued equivalence satisfying property (\ref{nv1}), similarly as in e.g., \cite{Bel2}, the difference is in the notion of  $\mu$-reflexivity.

 An additional important reason for our choice of $\mu$-reflexivity instead of the classical one is explained by Remark \ref{refl1} in Section \ref{omalg}.
\end{rem}

A \textbf{lattice-valued subalgebra} of an algebra $\mathcal{Q}=(Q,F)$ (here an \textbf{$L$-valued subalgebra} of $\mathcal{Q}$)  is a
function $\mu:Q\rightarrow L$ which is not constantly equal to 0, and which fulfils the following:
For any operation $f:Q^n\rightarrow Q$ from $F$ with arity $n > 0 \enskip (n \in \mathbb{N})$, 
and for all $a_1,\ldots ,a_n\in Q$, we have that
\begin{eqnarray} &\displaystyle\bigwedge_{i=1}^n\mu (a_i)\ls \mu (f(a_1,\ldots ,a_n)),&\label{compo}\\
&\mbox{and for a nullary operation } c\in F,
\;\mu (c)=1.&\label{cons}\end{eqnarray}

How the term operations behave in the lattice valued settings is formulated in the sequel. The proof goes easily by induction on the complexity of  terms.
\begin{prop}\label{termo} \sl Let $\mu:Q\rightarrow L$ be an $L$-valued subalgebra of an algebra $\mathcal{Q}$ and let $t(x_1,\ldots ,x_n)$ be a term in the language of $\mathcal{Q}$. If  $a_1,\ldots ,a_n\in Q$, then the following holds:
\begin{equation}\bigwedge_{i=1}^n\mu (a_i)\ls \mu (t(a_1,\ldots ,a_n)).\label{subalg}\end{equation} $\hspace*{\fill}\Box$\end{prop}

An $L$-valued relation $R :Q^2\rightarrow L$ on an algebra ${\mathcal{Q}}=(Q,F)$ is \textbf{compatible} with the operations in $F$ if the following two conditions holds:
for every $n$-ary operation $f\in F$, for all $a_1,\ldots,a_n, b_1,\ldots ,b_n\in Q$, and for every constant (nullary operation) $\;c \in F$
\begin{eqnarray}&&\bigwedge_{i=1}^n R (a_i,b_i)\ls R (f(a_1,\ldots,a_n),f( b_1,\ldots ,b_n));\label{comp1}\\&& R (c,c)=1.\label{comp2}\end{eqnarray}

\subsection{$L$-set}

The following is defined in \cite{FS} under the name of $\O$-set, and then adopted to a fuzzy framework in \cite{fugr}. In \cite{FS} $\O$ was a Heyting lattice, and $\O$-sets were used for modeling intuitionistic logic.

An \textbf{$L$-set} is a pair $(Q,E)$, where $Q$ is a nonempty set, and $E$  is  a symmetric and transitive  $L$-valued relation on $Q$, fulfilling the property (\ref{nv1}).

For an $L$-set $(Q,E)$, we denote by $\mu$ the $L$-valued function on $Q$, defined by
\begin{equation}\mu (x):=E(x,x).\label{mul}\end{equation}
We say that $\mu$ is \textsl{determined by $E$}.
Clearly, by the strictness property, $E$ is an $L$-valued relation on $\mu$, namely, it is an $L$-valued equality on $\mu$. That is why we say that in an $L$-set  $(Q,E)$, $E$ is an \textbf{$L$-valued equality} and $\mu (x)$ is the \textbf{degree of belonging} of $x$ to this $L$-set.

\begin{lemma}\sl\label{pomp} If $(Q,E)$ is an $L$-set and $p\in L$, then the cut $E_p$ is an equivalence relation on the corresponding cut $\mu_p$ of $\mu$.
\end{lemma}

\subsection{$L$-algebra; identities}\label{omalg}

Next we introduce a notion of a lattice-valued algebra with a lattice-valued equality.

Let  ${\mathcal{Q}}=(Q,F)$ be an algebra and $E:Q^2\rightarrow L$ an $L$-valued equality on $Q$, which is compatible with the operations in $F$. Then we say that $(\mathcal{Q},E)$ is an \textbf{$L$-algebra}.  Algebra ${\mathcal{Q}}$
 is the \textbf{underlying algebra} of $(\mathcal{Q},E)$.

 Now we present some cut properties of $L$-algebras.  These have been proved in \cite{fugr}, in the framework of groups.
\begin{prop}\label{nivoi}\sl Let $(\mathcal{Q},E)$ be an $L$-algebra. Then the following hold:

$(i\;)$ The function  $\mu:Q\rightarrow L$ determined by $E$ ($\mu (x)=E(x,x)$ for all $x\in Q$), is an $L$-valued subalgebra of $Q$.

 $(ii\;)$ For every $p\in L$, the cut $\mu_p$ of $\mu$ is a subalgebra of $\mathcal{Q}$, and

 $(iii\;)$ For every $p\in L$, the cut $E_p$ of $E$ is a congruence relation on $\mu_p$.
\end{prop}

Next we define how identities hold on  $L$-algebras, according to the \cite{BSATie}.

Let  $u(x_1,\ldots ,x_n)\approx v(x_1,\ldots ,x_n)$ (briefly $u\approx v$) be an identity in the type of an $L$-algebra $(\mathcal{Q},E)$. We assume, as usual, that variables appearing in terms $u$ and $v$ are from $x_1,\ldots ,x_n$  Then,  $(\mathcal{Q},E)$
\textbf{satisfies identity} $u\approx v$ (i.e., this
identity  \textbf{holds} on $(\mathcal{Q},E)$)  if the following condition is fulfilled:  %
\begin{equation}\bigwedge_{i=1}^n\mu (a_i)\ls
E(u(a_1,\ldots ,a_n),v(a_1,\ldots ,a_n)),\label{i1}\end{equation}
for all $a_1,\ldots, a_n\in Q$.

If $L$-algebra $(\mathcal{Q},E)$  satisfies an identity, then  this identity  need not hold  on $\mathcal{Q}$. On the other hand, if the underlying algebra fulfills an identity then also the corresponding $L$-algebra does.

\begin{prop}\label{uvod}\sl {\rm \cite{fugr}} If an identity $u\approx v$  holds on an algebra $\mathcal{Q}$, then it also holds on an $L$-algebra $(\mathcal{Q},E)$. \end{prop}

\begin{thm}\label{cuts}\sl {\rm \cite{fugr}} Let $(\mathcal{Q}, E)$ be an $L$-algebra, and $\mathcal{F}$ a set of identities in the language of $\mathcal{Q}$. Then, $(\mathcal{Q}, E)$ satisfies all identities in $\mathcal{F}$ if and only if for every $p\in L$
 the quotient algebra  $\mu_p{/}E_p$ satisfies the same identities.
\end{thm}

\begin{rem}\label{refl1}\rm
The fact that an $L$-algebra satisfy an identity while the same identity need not hold on the underlying algebra is caused by  $\mu$-reflexivity of the $L$-valued equality $E$. Namely, if $E$ would be reflexive in the classical sense ($E(x,x)=1$, for all $x\in Q$), then the cuts $\mu_1$ and $E_1$ would be the whole set $Q$ and the classical equality, respectively. Therefore, the quotient structure $\mu_1{/}E_1$ would be isomorphic to the underlying algebra $\mathcal{Q}$. By Theorem \ref{cuts}, in this case $L$-algebra $(\mathcal{Q},E)$ would not bring anything new, it would simply repeat properties satisfied by algebra $\mathcal{Q}$.
\end{rem}
\section{Results}
\subsection{$L$-groupoid, $L$-quasigroup}
Let $L$ be a complete lattice. According to the definition of an $L$-algebra, an \textbf{$L$-groupoid} is a structure $(\mathcal{Q},E)$, where $\mathcal{Q}=(Q,\cdot )$ is a groupoid and $E:Q^2\rightarrow L$ an $L$-valued compatible equality over  $\mathcal{Q}$.

Let $(\mathcal{Q},E)$ be an $L$-groupoid. Each of the formulas $a\cdot x = b$ and $y\cdot a = b$,  $a,b\in Q$, $x,y$ -- variables, is  a \textbf{linear equation over $(\mathcal{Q},E)$}.

We say that an equation $a\cdot x = b$ is \textbf{solvable over $(\mathcal{Q},E)$} if there is $c\in Q$ such that
\begin{equation}
\mu (a)\wedge\mu (b)\ls\mu (c)\wedge E(a\cdot c, b).\label{jea}\end{equation}
Analogously, an equation  $y\cdot a = b$  is \textbf{solvable over $(\mathcal{Q},E)$} if there is $d\in Q$ such that
\begin{equation}
\mu (a)\wedge\mu (b)\ls\mu (d)\wedge E(d\cdot a, b).\label{jeb}\end{equation}

Elements $c$ and $d$ are \textbf{solutions of equations $a\cdot x = b\;$ and $\;y\cdot a = b$, respectively in $(\mathcal{Q},E)$}.

If $c$ and $d$ are solutions of $a\cdot x = b\;$ and $\;y\cdot a = b$, respectively in $(\mathcal{Q},E)$, then obviously  for every $p\in L$ satisfying $p\ls\mu (a)\wedge\mu (b)$, we have
\begin{equation}
p\ls\mu (c)\wedge E(a\cdot c, b),\label{jea'}\end{equation}
and
\begin{equation}
p\ls\mu (d)\wedge E(d\cdot a, b).\label{jeb'}\end{equation}

Each of the above equations is \textbf{$E$-uniquely solvable over $(\mathcal{Q},E)$} if  the following hold:

If $c$ is  a solution of the equation $a\cdot x = b$ over $(\mathcal{Q},E)$ and   $c_1\in Q$ fulfills $E(a\cdot c_1,b)\gs p$ for some $p\ls\mu (a)\wedge\mu (b)$, then
\begin{equation}E(c,c_1)\gs p.\label{dva}\end{equation}
Analogously, if $d$ is  a solution of the equation $y\cdot a = b$ over $(\mathcal{Q},E)$ and   $d_1\in Q$ fulfills $E(d_1\cdot a,b)\gs p$ for some $p\ls\mu (a)\wedge\mu (b)$, then
\begin{equation}E(d,d_1)\gs p.\label{dvb}\end{equation}

  If $c_1$ and $d_1$ are (additional) solutions of equations $a\cdot x=b$ and $y\cdot a=b$, respectively, then clearly conditions (\ref{dva}) and (\ref{dvb}) hold. Hence,
 an $E$-uniquely solvable equation may have several solutions. All these solutions are equal up to the $ L$-equality $E$. More precisely, we have the following.

\begin{thm}\label{prt}\sl
Let $(\mathcal{Q},E)$ be an $L$-groupoid. If equations $a\cdot x = b$ and $y\cdot a = b$, are $E$-uniquely solvable over $(\mathcal{Q},E)$ for all $a,b\in Q$, then for every $p\in L$ the quotient groupoid $\mu_p{/}E_p$ is a quasigroup.
\end{thm}
\proof
Let $p\in L$, and let $a,b\in\mu_p$. Then obviously $p\ls\mu (a)\wedge\mu (b)$. Consider the equation $a\cdot x = b$. Then, by assumption, there is $c\in\mu_p$, such that  condition (\ref{jea}) is valid, and if $c_1\in Q$ fulfills $E(a\cdot c_1,b)\gs p$ for some $p\ls\mu (a)\wedge\mu (b)$ then (\ref{dva}) holds.
 By (\ref{jea}), also $E(a\cdot c, b)\gs p$, i.e., $(a\cdot c,b)\in E_p$. Since $E_p$ is a congruence over the subgroupoid $\mu_p$ of $\mathcal{Q}$, we get $[a\cdot c]_{E_p}=[b]_{E_p}$, i.e., $[a]_{E_p}\cdot [c]_{E_p}=[b]_{E_p}$. Therefore, an equation of the form $A\cdot X=B$, $A ,B\in\mu_p{/}E_p$ is solvable. By (\ref{dva}), the solution is unique in the classical sense. Indeed, if also $[a]_{E_p}\cdot [c_1]_{E_p}=[b]_{E_p}$, for some $c_1\in\mu_p$, then $[a\cdot c_1]_{E_p}=[b]_{E_p}$, and hence
$ p\ls E(a\cdot c_1,b).$ Therefore by (\ref{dva}), $E(c,c_1)\gs p$, hence $[c]_{E_p}=[c_1]_{E_p}$ and the solution is unique.

The proof that every equation of the form $Y\cdot A=B$ is also uniquely solvable over $\mu_p{/}E_p$ is analogous.
\endproof

We say that an $ L$-groupoid  $(\mathcal{Q},E)$ is an \textbf{$ L$-quasigroup}, if every equation of the form $a\cdot x = b$ or $y\cdot a = b$ is $E$-uniquely solvable over $(\mathcal{Q},E)$.

The converse of Theorem \ref{prt} also holds, as follows.

\begin{thm}\label{novo}\sl Let $(\mathcal{Q},E)$ be an $ L$-groupoid. If for all $a,b\in Q$ and for  every $p\ls\mu (a)\wedge\mu (b)$ the quotient groupoid $\mu_p{/}E_p$ is a quasigroup, then $(\mathcal{Q},E)$ is an $ L$-quasigroup.
\end{thm}
\proof
Let $a,b\in Q$ and let $p=\mu (a)\wedge\mu (b)$. By assumption, $\mu_p{/}E_p$ is a quasigroup, hence  the equations $[a]_{E_p}\cdot X=[b]_{E_p}$ and $Y\cdot [a]_{E_p}=[b]_{E_p}$ have unique solutions, $X=[c]_{E_p}$ and $Y=[d]_{E_p}$, for some $c,d\in\mu_p$. Hence,
$[a]_{E_p}\cdot [c]_{E_p}=[b]_{E_p}$ and $[d]_{E_p}\cdot [a]_{E_p}=[b]_{E_p}$, i.e., $[a\cdot c]_{E_p}=[b]_{E_p}$ and $[d\cdot a]_{E_p}=[b]_{E_p}$. Then $E(a\cdot c,b)\gs p$, and since $c\in\mu_p$,
\[p=\mu (a)\wedge\mu (b)\ls\mu (c)\wedge E(a\cdot c,b),\]
equation $a\cdot x=b$ is solvable over $(\mathcal{Q},E)$; similarly, equation $y\cdot a = b$ is also solvable over $(\mathcal{Q},E)$. These equations are $E$-uniquely solvable. Indeed, if $c$ is a solution of $a\cdot x=b$ over $(\mathcal{Q},E)$, and there is $c_1\in Q$ such that
 $E(a\cdot c_1,b)\gs q$ for some $q\ls p=\mu (a)\wedge\mu (b)$, then
$[a]_{E_q}\cdot [c_1]_{E_q}=[b]_{E_q}$. Since $q\ls p$, we have $\mu_p\subseteq\mu_q$ and $c\in\mu_q$. By assumption,
$\mu_q{/}E_q$ is a quasigroup, therefore $[a]_{E_q}\cdot [c]_{E_q}=[b]_{E_q}$ and $[c]_q=[c_1]_q$. Therefore, $E(c,c_1)\gs q$, and the equation $a\cdot x=b$ is $E$-uniquely solvable over $(\mathcal{Q},E)$. Analogously,
 the equation $y\cdot a=b$ is $E$-uniquely solvable over $(\mathcal{Q},E)$.

Therefore, $(\mathcal{Q},E)$ is an $ L$-quasigroup.
\endproof

\subsection{$ L$-equasigroup}
Let $\mathcal{Q}=(Q,\cdot, \backslash,/)$ be an algebra in the language with three binary operations, $ L$ a complete lattice and $E:Q^2\rightarrow L$ an $ L$-valued compatible equality over  $\mathcal{Q}$. Then, $(\mathcal{Q},E)$ is an \textbf{$ L$-equasigroup}, if identities $Q1,\ldots , Q4$ hold.
By (\ref{i1}), this means that the following formulas should be satisfied, where, as before, $\mu:Q\rightarrow L$ is defined by $\mu (x)=E(x,x)$:\\

$QE1:\;\mu (x)\wedge\mu (y)\ls E(y, x \cdot (x \backslash y)) ;$

$QE2:\;\mu (x)\wedge\mu (y)\ls E(y , x \backslash (x \cdot y)) ;$

$QE3:\;\mu (x)\wedge\mu (y)\ls E(y , (y / x) \cdot x) ;$

$QE4:\;\mu (x)\wedge\mu (y)\ls E(y , (y \cdot x) / x) .$

\begin{thm}\sl\label{oalg}
If $((Q,\cdot, \backslash,/),E)$ is an {$ L$-equasigroup}, then for every $p\in L$, the quotient structure $\mu_p{/}E_p$ is a classical equasigroup.
\end{thm}
\proof
This is a straightforward consequence of Theorem \ref{cuts}.
\endproof

\begin{cor}\sl
If $((Q,\cdot, \backslash,/),E)$ is an {$ L$-equasigroup}, then  $((Q,\cdot),E)$ is an $ L$-quasigroup.
\end{cor}
\proof
If $((Q,\cdot, \backslash,/),E)$ is an $ L$-equasigroup, then $((Q,\cdot,E)$ is an $ L$-groupoid, that is $E$ is an $ L$-equality on the groupoid $(Q,\,\cdot\,)$. Indeed, the operation $\,\cdot\,$ on $Q$ is the same one from  $((Q,\cdot, \backslash,/),E)$, hence compatibility of $E$ holds.

Next, by Theorem \ref{oalg} every structure $\mu_p{/}E_p$ is a classical equasigroup. By Theorem \ref{equiv} every such structure is a quasigroup, hence by Theorem \ref{prt}, $ L$-groupoid $((Q,\cdot),E)$ is an $ L$-quasigroup.
\endproof

The converse follows by the Axiom of Choice (AC).

 Let $((Q,\cdot),E)$ be an $ L$-groupoid which is an $ L$-quasigroup.
 By Theorem \ref{prt}, for every $p\in L$, the quotient groupoid $(\mu_p{/}E_p, \cdot )$ is a quasigroup, where the operation $\cdot $ is  defined by  $[a]_{E_p}\cdot [b]_{E_p}=[a\cdot b]_{E_p}$,  $a,b\in\mu_p.$
By Theorem \ref{equiv}, the structure  $(\mu_p{/}E_p, \,\cdot\,,\backslash\,,/\,)$ is an equasigroup, where the operations $\,\backslash\,$ and $\,/\,$ are the usual ones:
\begin{eqnarray*} &&[a]_{E_p}\backslash\, [b]_{E_p}=[c]_{E_p}\;\mbox{ if and only if }\; [a]_{E_p}\cdot [c]_{E_p} =[b]_{E_p}, \;\mbox{ and }\\ &&[b]_{E_p}/\,[a]_{E_p}=[d]_{E_p}\;\mbox{ if and only if }\; [d]_{E_p}\cdot [a]_{E_p} =[b]_{E_p}.\end{eqnarray*}
Let us define binary operations  $ \backslash$ and $/$ over $Q$ in the following way:

 For every pair $a,b\in Q$,
$\; a\backslash\,b=c,\;$ where $c$ is an element chosen by AC from $[a]_{E_p}\backslash\, [b]_{E_p}$ in the quasigroup $\mu_p{/}E_p$, where $p=\mu (a)\wedge\mu (b)$. Analogously, $\; b/\,a=d,\;$ where $d$ is  chosen by the AC from $[b]_{E_p}/\, [a]_{E_p}$ in $\mu_p{/}E_p$, for $p=\mu (a)\wedge\mu (b)$.

\begin{lemma}\sl
 Let $((Q,\cdot),E)$ be an $ L$-groupoid which is an $ L$-quasigroup. Then
 the operations $ \backslash$ and $/$ over $Q$ are well defined.
\end{lemma}
\proof
Let  $a,b\in Q$,
$\; a\backslash\,b=c,\;$ where $c$ is an element chosen by AC from $[a]_{E_p}\backslash\, [b]_{E_p}$ in the quasigroup $\mu_p{/}E_p$, where $p=\mu (a)\wedge\mu (b)$.  Elements $a$ and $b$ belong to $\mu_p$, since $a\in\mu_p$ if and only if $\mu (a)\gs p=\mu (a)\wedge\mu (b)$, and the latter obviously holds, similarly for $b$. Since $(\mu_p{/}E_p, \,\cdot\,,\backslash\,,/\,)$ is an equasigroup, the class $[a]_{E_p}\backslash\, [b]_{E_p}$ exists. Therefore, there is $c\in\mu_p$, and being a chosen element, it is unique. Similarly, one can show that also the operation $/\,$ is well defined.
\endproof

\begin{lemma}\label{lnov}\sl
 Let $((Q,\cdot),E)$ be an $ L$-groupoid which is an $ L$-quasigroup. Then for every $q\in L$ and for all $a,b\in\mu_q$,  in the quasigroup $(\mu_q{/}E_q, \cdot\,,\backslash\,,{/}\,)$ we have
 $ [a\backslash b]_{E_q}=[a]_{E_q}\backslash [b]_{E_q},
 $ and
 $[a/ b]_{E_q}=[a]_{E_q}/ [b]_{E_q},
 $
 where the operations $\backslash$ and $/$ on the left hand sides are the ones defined on $Q$ by AC.
\end{lemma}
\proof
Let $a,b\in\mu_q$, and $p=\mu (a)\wedge\mu (b)$. Then, in the quasigroup $(\mu_p{/}E_p, \cdot\,,\backslash\,,{/}\,)$ we have

$[a]_{E_p}\cdot [a\backslash b]_{E_p}=[b]_{E_p}$, i.e.,
$[a\cdot (a\backslash b)]_{E_p}=[b]_{E_p}$.

\noindent Since $q\ls p$, we have $\mu_p\subseteq\mu_q$ and $E_p\subseteq E_q$, thus we get

$[a\cdot (a\backslash b)]_{E_q}=[b]_{E_q}$.

$E_q$ is a congruence relation on the subgroupoid $\mu_q$,  therefore it is compatible with the operation $\,\cdot\,$, implying

$[a]_{E_q}\cdot [a\backslash b]_{E_q}=[b]_{E_q}$ on the quasigroup $\mu_q{/}E_q$.

In the quasigroup $\mu_q{/}E_q$, the class $[a\backslash b]_{E_q}$ is the unique satisfying the above equality, for given $a,b$. Moreover, this class is precisely the one obtained as a result of the application of $\backslash$ on $[a]_{E_q}$ and $[b]_{E_q}$:

$[a]_{E_q}\backslash [b]_{E_q}=[a\backslash b]_{E_q}$.

The proof for the remaining operation $/$ is analogous.
\endproof

\begin{thm}\sl Let $((Q,\cdot),E)$ be an $ L$-groupoid which is an $ L$-quasigroup. Then   the structure $((Q,\,\cdot\,,\backslash\,,{/}\,),E)$ is an $ L$-equasigroup, where the binary operations $ \backslash$ and $/$ over $Q$ are defined by Axiom of Choice as above.
\end{thm}
\proof
Suppose that $\mathcal{Q}=(Q,\cdot)$ is a groupoid and $(\mathcal{Q},E)$  an $ L$-quasigroup, with $E:Q^2\rightarrow L$ being a compatible $ L$-valued equality on $\mathcal{Q}$.

We prove that $((Q,\,\cdot\,,\backslash\,,{/}\,),E)$ is an $ L$-algebra, moreover that it is an $ L$-equasigroup, where $\,\cdot\,$ is the starting operation in the groupoid $(Q,\,\cdot\,)$ while $\,\backslash\,$ and $\,/\,$ are operations defined above by the use of the Axiom of Choice.

To prove that this structure is an $ L$-algebra, we have to show that $E$ is compatible with new operations $\,\backslash\,$ and $\,/\,$ (it is already compatible with  $\,\cdot\,$). Indeed, for every $q\in L$, by Proposition \ref {nivoi} the cut $E_q$ is an ordinary congruence on the groupoid $(\mu_q,\,\cdot\,)$. In addition, the restrictions of the new binary operations $\,\backslash\,$ and $\,/\,$ to $\mu_q$ are also operations on this set: If $a,b\in\mu_q$, then $\mu (a)\gs q$ and $\mu (b)\gs q$.
Hence $\mu (a)\wedge\mu (b)\gs q$ and therefore $\mu_p\subseteq\mu_q$, for $p=\mu (a)\wedge\mu (b)$. Obviously, $a,b\in\mu_p$, and by the definition of the new operations we have also that $a\backslash \,b\in\mu_p$  and $a/\,b\in\mu_p$. Since $\mu_p$ is a subset of $\mu_q$, it follows that $a\backslash\, b\in\mu_q$  and $a/\,b\in\mu_q$, proving that these restrictions are operations on $\mu_q$. Consequently, we have an algebra $(\mu_q,\,\cdot\,,\backslash\,,\,/\,)$ and $E_q$ is an equivalence relation on it, compatible with the first of these three binary operations. Compatibility with remaining two:   if $x,y,u,v\in\mu_q$ and $(x,y),(u,v)\in E_q$, then
 by Lemma \ref{lnov}, $x\backslash\, u\in [x]_{E_q}\backslash\, [u]_{E_q}$, i.e.,
\[[x]_{E_q}\backslash\, [u]_{E_q}=[x\backslash\, u]_{E_q}
\;\;\mbox{ and similarly }\;\;[y]_{E_q}\backslash\, [v]_{E_q}=[y\backslash\, v]_{E_q}.\]
 But since $E_q$ is an equivalence relation on $\mu_q$, we have also $[x]_{E_q}= [y]_{E_q}$ and $[u]_{E_q}= [v]_{E_q}$. Therefore
\[[x\backslash\, u]_{E_q} =[y\backslash\, v]_{E_q}\;\;\; \mbox{ and }\;\;\;(x\backslash\, u,y\backslash\, v)\in E_q.\]
Analogously, one could prove that $E_q$ is compatible with the restriction of the operation $\,/\,$ to $\mu_q$. Hence, for every $q\in L$, the cut $E_q$ is a congruence on $\mu_q$, hence $E$ is an $ L$-valued equality on the algebra $(Q,\,\cdot\,,\backslash\,,{/}\,)$. In this way we have proved that  $((Q,\,\cdot\,,\backslash\,,{/}\,),E)$ is an $ L$-algebra. Finally, we prove that it is an $ L$-equasigroup. For every $q\in L$, $(\mu_p{/}E_p, \,\cdot\,,\backslash\,,/\,)$ is an equasigroup, satisfying identities $Q1$ -- $Q4$. By  Theorem \ref{cuts}, the corresponding $ L$-algebra $((Q,\,\cdot\,,\backslash\,,{/}\,),E)$ also satisfies these identities, i.e., formulas $QE1$ -- $QE4$ hold. Therefore, $((Q,\,\cdot\,,\backslash\,,{/}\,),E)$ is an $ L$-equasigroup.
\endproof
\subsection{Example}
We present a toy example in which a groupoid equipped with a fuzzy equality is an $L$-quasigroup. By this example we also illustrate the procedure of solving linear equation w.r.t. fuzzy equality.

Let $(Q,\,\cdot\,)$ be a groupoid given in Table 1. Obviously, this groupoid  is not a quasigroup, e.g., equation $a\cdot x=d$, as visible from the table, does not have a solution in $Q$.
\begin{center}
\textit{\begin{tabular}{c|ccccc}$\cdot$&a&b&c&d&e \\ \hline a&b&c&a&a&e\\b&a&b&c&d&e\\c&c&a&b&b&e\\d&d&a&b&b&e\\e&e&e&e&e&a\end{tabular}
}
\end{center}
\begin{center} \textbf{Table 1} \end{center}
The lattice $ L$ is given by the diagram in Figure 1, and an $ L$-valued equality is presented by Table 2. Hence, $((Q,\,\cdot\,),E)$ is an $ L$-groupoid.

\begin{center}
\unitlength 1mm 
\linethickness{0.4pt}
\ifx\plotpoint\undefined\newsavebox{\plotpoint}\fi 
\begin{picture}(49,63)(0,-5)
\put(22,55){\circle*{2}}
\put(22,55){\line(-1,-1){12}}
\put(10,43){\line(1,-1){12}}
\put(22,31){\line(1,1){12}}
\put(34,43){\line(-1,1){12}}
\put(10,43){\circle*{2}}
\put(22,31){\circle*{2}}
\put(34,43){\circle*{2}}
\put(34,43){\line(1,-1){12}}
\put(46,31){\line(0,-1){10}}
\put(46,31){\circle*{2}}
\put(46,21){\circle*{2}}
\put(46,21){\line(-1,-1){12}}
\put(34,9){\line(-1,1){12}}
\put(22,21){\line(0,1){10}}
\put(22,21){\circle*{2}}
\put(34,9){\circle*{2}}
\put(22,58){\makebox(0,0)[cc]{$1$}}
\put(7,43){\makebox(0,0)[cc]{$q$}}
\put(37,44){\makebox(0,0)[cc]{$p$}}
\put(19,29){\makebox(0,0)[cc]{$r$}}
\put(49,31){\makebox(0,0)[cc]{$u$}}
\put(19,20){\makebox(0,0)[cc]{$w$}}
\put(49,20){\makebox(0,0)[cc]{$v$}}
\put(34.3,5.9){\makebox(0,0)[cc]{$0$}}
\put(17,7){\makebox(0,0)[cc]{Lattice $ L$}}
\put(28,0){\makebox(0,0)[cc]{\textbf{Figure 1}}}
\end{picture}
\end{center}

\begin{center}
\textit{\begin{tabular}{c|ccccc}E&a&b&c&d&e \\ \hline a&$1$&p&p&r&v\\b&p&$1$&p&r&v\\c&p&p&$1$&q&v\\d&r&r&q&q&$0$\\e&v&v&v&$0$&u\end{tabular}
}
\end{center}
\begin{center} \textbf{Table 2} \end{center}

The function  $\mu:Q\rightarrow L$  ($\mu (x)=E(x,x)$ for all $x\in Q$):
\begin{center}
$\mu=\left(\begin{array}{ccccc}a&b&c&d&e\\1&1&1&q&u\end{array}\right).$
\end{center}
The subgroupoids of $((Q,\,\cdot\,),E)$, which are cuts of $\mu$:

$\mu_1=\mu_p=\{ a,b,c\}$,

$\mu_q=\mu_r=\mu_w=\{ a,b,c,d\}$,

$\mu_u=\mu_v=\{ a,b,c,e\}$,

$\mu_0=\{ a,b,c,d,e\}$.

\noindent The quotient groupoids over the corresponding cuts of $E$ are the following:

$\mu_1{/}E_1=\{\{ a\},\{ b\},\{ c\}\}$,

 $\mu_p{/}E_p=\{\{ a, b, c\}\}$,

$\mu_q{/}E_q=\{\{ a\},\{ b\},\{ c,d\}\}$,

$\mu_r{/}E_r=\mu_w{/}E_w=\{\{ a, b, c,d\}\}$,

 $\mu_u{/}E_u=\{\{ a,b , c\},\{ e\}\}$,

  $\mu_v{/}E_v=\{\{ a,b , c, e\}\}$,

  $\mu_0{/}E_0=\{\{ a,b , c, d,e\}\}$.

All these quotient structures are quasigroups, hence the starting $L$-groupoid is an $L$-quasigroup, and every linear equation is $E$-uniquely solvable over it.
E.g., the mentioned equation $a\cdot x=d$  which does not have a classical solution in $Q$, possesses a solution with respect to fuzzy equality $E$. Indeed, due to $\mu (a)\wedge\mu (d)=q$, this solution is element $b$, since the class $X=\{ b\}$ is the unique solution of the equation $[a]_{E_q}\cdot X=[d]_{E_q}$ over the quasigroup $\mu_q{/}E_q$ (observe that $[d]_{E_q}=\{ c,d\}$). By (\ref{jea}), we have
\[\mu (a)\wedge\mu (d)=q\ls \mu (b)\wedge E(a\cdot b,d)= \mu (b)\wedge E(c,d)=1\wedge q=q.\]
Hence, $a\cdot b$ and $d$ are $E$-equal with grade $q$.

\subsection{$L$-loop and $L$-group}
As defined in \cite{fugr}, an $L$-algebra $(\mathcal{G},E)$ is an \textbf{$L$-group}, if the underlying algebra $\mathcal{G}=(G,\cdot , {\;}^{-1},e)$ has a binary operation $\,\cdot\,$, a unary operation ${\;}^{-1}$, a constant $e$, and  the following formulas hold:

$LG1:\;$ $ \mu (x)\wedge\mu (y)\wedge\mu (z)\ls E(x\cdot (y\cdot z), (x\cdot y)\cdot z);$

 $LG2:\;$ $\mu (x)\ls E(x\cdot e, x),\;\; \mu (x)\ls E(e\cdot x, x);$

 $LG3:\;$ $\mu (x)\ls E(x\cdot x^{-1}, e),\;\; \mu (x)\ls E( x^{-1}\cdot x, e).$

The following is a consequence of Theorem \ref{cuts}.
\begin{thm}\label{lgr}\sl An $L$-algebra $((G,\cdot , {\;}^{-1},e),E)$ is an $L$-group if and only if for every $p\in L$, the quotient cut-subalgebra $\mu_p{/}E_p$ is a group.
\end{thm}

Observe that $e$ corresponds to the constant in the language, therefore $E(e,e)=\mu (e)=1$.
Now, if $(\mathcal{G},E)$ is an $L$-group, then by the condition (\ref{nv1}), we get $E(e,x)< 1$ whenever $x\neq e$, and thus by $LG2$,
\[ 1=E(e,e)\ls E(e\cdot e,e).\]
Hence, in the underlying algebra $\mathcal{G}$, $e\cdot e=e$.

We define an \textbf{$L$-loop} as an $L$-algebra $(\mathcal{Q},E)$, where $\mathcal{Q}=(Q,\cdot,e)$ is a structure with a binary operation $\,\cdot\,$ and a constant $e$, $((Q,\,\cdot\,),E)$ is an $L$-quasigroup,  $E(e,e)=1$ and the formula $LG2$ holds. %

An \textbf{$L$-semigroup} \cite{fi1} is an $L$-algebra $((Q,\,\cdot\,),E)$ where $(Q,\,\cdot\,)$ is a groupoid and the formula $LG1$ holds.

The proof of the following theorem depends on the Axiom of Choice (AC).
\begin{thm}\sl
Let $((Q,\cdot,e),E)$ be an $L$-algebra. There is a unary operation ${\;}^{-1}$ on $Q$ such that $((Q,\cdot , {\;}^{-1},e),E)$ is an $L$-group if and only if $((Q,\cdot),E)$ is an $L$-semigroup and $((Q,\cdot ,e),E)$  an $L$-loop.
\end{thm}
\proof
Let $((Q,\cdot,e),E)$ be an $L$-algebra and suppose there is a unary operation ${\;}^{-1}$ on $Q$ such that $((Q,\cdot , {\;}^{-1},e),E)$ is an $L$-group. Then by Theorem \ref{lgr}, for every $p\in L$ $\mu_p{/}E_p$ is a group, hence it is a semigroup and a quasigroup.
Then clearly, $((Q,\cdot),E)$ is an $L$-semigroup by Theorem \ref{cuts}  and an $L$-quasigroup by Theorem \ref{novo}. By $LG2$, $((Q,\cdot ,e),E)$  is an $L$-loop.
%








To prove the converse, we assume that $((Q,\cdot,e),E)$ is an $L$-algebra, such that $((Q,\cdot),E)$ is an $L$-semigroup and $((Q,\cdot ,e),E)$  an $L$-loop. We define a unary operation ${\;}^{-1}$ on $Q$ as follows. Let $a\in Q$, such that $\mu (a)=p$. By assumption and by Theorem \ref{cuts}, $\mu_p{/}E_p$ is a loop with the identity element $[e]_{E_p}$. Therefore, equation $[a]_{E_p}\cdot X=[e]_{E_p}$ has a unique solution in $\mu_p{/}E_p$, the class $[c]_{E_p}$, for some $c\in\mu_p$. Now, by AC we define $a^{-1}$ to be an arbitrary element in $[c]_{E_p}$. This operation is well defined, since for every $a\in Q$,  equation $[a]_{E_p}\cdot X=[e]_{E_p}$ has a unique solution, a class in $\mu_p{/}E_p$, for $p=\mu (a)$; the chosen element from the corresponding class is unique by construction. In addition, since $E_p$ is a congruence on $\mu_p$, $[a]_{E_p}\cdot [a^{-1}]_{E_p}=[e]_{E_p}$ and $p=\mu (a)$, we have
\[\mu (a)\ls E(a\cdot a^{-1},e), \;\;\mu (a)\ls E(a^{-1}\cdot a,e),\]
proving $LG3$. $LG1$ and $LG2$ hold by assumption, hence  $((Q,\cdot , {\;}^{-1},e),E)$ is an $L$-group.
\endproof

\begin{rem}\rm
Let us mention that the equivalence among $L$-groups and associative $L$-loops essentially depends on the language in which these structures are defined. An option, like in the classical algebra, could be that the underlying structures were groupoids without a nullary operation in the language (identity element being required to exist in the groupoid). However, in the framework of $L$-algebras, the quotient structures over cuts would not necessarily share the same identity element, and the equivalence would not be fulfilled. Examples are easy to construct. E.g., it could be any groupoid having two disjoint subgroupoids which are groups. With a suitable $L$-valued equality, these subgroupoids could be cuts, and the quotient structures would become disjoint groups. Hence, no common identity could exist.  
\end{rem}

\subsection{An application}

As presented above, in order to be able to find unique solutions of linear equations, we do not need a quasigroup, only a groupoid is needed. Quasigroups then appear as quotients over cuts with respect to a fuzzy equality. This is much weaker requirement for the starting binary operation. We can go further, relaxing also this weaker requirement of quasigroups over quotients.
A motivation comes from applications, as follows.

Namely, let $\,\cdot\,$ be an arbitrary binary operation on a set $Q$. In financial transactions, in managing data etc., it is frequently necessary to solve the equation $a\cdot x=b$ for particular (not necessarily all) $a,b\in Q$. In real situations it may happen that the solution does not exists, or it might be impossible to identify strictly equal objects. In such situations we do not deal with classical equality (\,=\,), but elements, objects in $Q$ might be equal 'up to some extent', in which case we have an appropriate fuzzy equality $E$. Usually, such an equality, i.e., its membership values are known, calculated in advance, depending on the context. So, we may wish to find "fuzzy" solutions, i.e., element(s) $c$ for which, intuitively, \[\mbox{ $a\cdot c\;$ equals $\,b\;$ with respect to the fuzzy equality $E$.}\]

According to the above explanation, we introduce the following definition.

Let   $\mathcal{Q}=(Q,\,\cdot\,)$ be an arbitrary groupoid and $E: Q^2\rightarrow L$ an $ L$-valued equality over $\mathcal{Q}$;  let also $a,b\in Q$. Then we say that the equation $a\cdot x=b$ has a \textbf{unique solution w.r.t.\ $E$}, if this equation is $E$-uniquely solvable over the $ L$-groupoid $(\mathcal{Q},E)$.

The following theorem is not a direct consequence of Theorem \ref{novo}, still the proof uses the same arguments. 
\begin{thm}\label{posl}\sl Let $\mathcal{Q}=(Q,\,\cdot\,)$ be an arbitrary groupoid, let $a,b$ be particular elements in $Q$, and let $E:Q^2\rightarrow L$ be an $ L$-valued equality over $\mathcal{Q}$.  Then the equation  $a\cdot x=b$ has a unique solution w.r.t.\ $E$, if the equation $[a]_{E_p}\cdot X=[b]_{E_p}$, for $p=\mu (a)\wedge\mu (b)$, has  a (classical) unique solution in the quotient groupoid $(\mu_p{/}E_p, \,\cdot\,)$.
\end{thm}

What does unique solvability in this context practically means? As mentioned, in real situations, for chosen $a,b$ in $Q$ the equation $a\cdot x=b$ may not have any (classical) solution. Still, there might be some "close" values, with respect to the fuzzy equality $E$. Theorem \ref{posl} tells us that the most close solutions w.r.t.\ fuzzy equality $E$ are elements of the class $C\in\mu_p{/}E_p$ , $p=\mu (a)\wedge\mu (b)$, which is a classical solution of the equation $[a]_{E_p}\cdot X=[a]_{E_p}$. In this case we consider every element $c\in C$ to be a solution of the equation $a\cdot x=b$ in the groupoid $\mathcal{Q}$.

\section{Conclusion}
This investigation is focussed to classical linear equations with one operation, appearing frequently in real problems. In solution procedures we use fuzzy (lattice valued) equality and cut techniques. The background for our research is the general algebra and so-called $L$-quasigroups, a generalization of the classical structures in which these equations have unique solutions. $L$-quasigroups are equipped with a fuzzy equality, with respect to the basic operation they are not quasigroups, hence being much closer to structures appearing in real applications.

Our technique is new and widely applicable. Developing this procedure it could be possible to deal with linear equations with two operations and several unknowns. Consequently, we intend to focuss on classical systems of linear equations, in the situations where not all data  are known and the classical equality has to be replaced by a fuzzy one.

\end{document}